\documentclass{article}
\usepackage[dvipdfmx]{graphicx}
\usepackage{spconf}
\usepackage{latexsym}
\usepackage{amsmath}
\usepackage{amssymb}
\usepackage{amsfonts}
\usepackage{amsthm}
\usepackage{setspace}
\usepackage[ruled,vlined]{algorithm2e}
\usepackage{url}
\usepackage{cite}
\usepackage{balance}

\makeatother

\theoremstyle{plain}
\newtheorem{prop}{Proposition}

\newtheorem{remk}{Remark}

\newcommand{\argmin}{\mathop{\rm argmin}\limits}

\graphicspath{{./images/}}
\def\x{{\mathbf x}}
\def\xii{{\boldsymbol \xi}}

\def\y{{\mathbf y}}

\def\w{{\mathbf w}}

\def\z{{\mathbf z}}
\def\v{{\mathbf v}}

\def\e{{\mathbf e}}

\def\n{{\mathbf n}}
\def\t{{\mathbf t}}
\def\u{{\mathbf u}}

\def\bzeta{{\boldsymbol \zeta}}

\def\bvarepsilon{{\boldsymbol \varepsilon}}
\def\A{{\mathbf A}}

\def\D{{\mathbf D}}

\def\R{{\mathbb R}}

\def\I{{\mathbf I}}

\def\PHI{{\mathbf \Phi}}
\def\PSI{{\mathbf \Psi}}

\def\LAMBDA{{\mathbf \Lambda}}

\def\Rmath{{\mathcal R}}

\def\Sfrak{{\mathfrak S}}

\def\Scal{{\mathcal S}}
\def\Vcal{{\mathcal V}}

\DeclareMathOperator{\prox}{prox}

\title{Efficient Constrained Signal Reconstruction\\ by Randomized Epigraphical Projection}

\name{Shunsuke Ono
\thanks{
The work was partially supported by JSPS Grants-in-Aid (17K12710) and JST-PRESTO.}}
\address{
Tokyo Institute of Technology
}

\begin{document}
\maketitle

\begin{abstract}
This paper proposes a randomized optimization framework for constrained signal reconstruction, where the word ``constrained'' implies that data-fidelity is imposed as a hard constraint instead of adding a data-fidelity term to an objective function to be minimized.
Such formulation facilitates the selection of regularization terms and hyperparameters, but due to the non-separability of the data-fidelity constraint, it does not suit block-coordinate-wise randomization as is.
To resolve this, we give another expression of the data-fidelity constraint via epigraphs, which enables to design a randomized solver based on a stochastic proximal algorithm with randomized epigraphical projection.
Our method is very efficient especially when the problem involves non-structured large matrices. 
We apply our method to CT image reconstruction, 
where the advantage of our method over the deterministic counterpart is demonstrated.
\end{abstract}
\begin{keywords} signal reconstruction, constrained optimization, stochastic optimization, epigraphical projection
\end{keywords}

\section{Introduction}\label{intro}
Signal reconstruction from incomplete and/or degraded observation is a fundamental problem arising from various applications, such as medical imaging, microscopy, tomography, spectral imaging and computational photography.
Such a problem is often reduced to a convex optimization problem that involves a regularization term, modeling some desirable properties on the signal of interest, and a data fidelity term, enforcing consistency with observed data.

Proximal splitting algorithms \cite{TechCombettes,Boydprox} have played a central role in solving such problems.
These methods, especially primal-dual splitting type algorithms \cite{PDChambolle,PDCombettes,PDCondat,PDSVu,PDSsurvey,SPL2017}, are efficient in the sense that they require only simple operations like matrix-vector multiplications and evaluation of proximity operators.
However, even such a simple operation becomes computationally expensive in many applications.
A typical case is computed tomography (CT), where a matrix representing the observation process is large and is not structured\footnote{Here the word ``structured'' means that the matrix-vector multiplication can be computed efficiently via some operation. An instance of structured matrices is a uniform blur matrix, which can be diagonalized by FFT.}, resulting in large computational costs and memory requirements \cite{CTTV,CTreview}.

Recently, \emph{stochastic primal-dual splitting algorithms} have been intensively studied for stochastic optimization \cite{SPDSPesquet,SPDSZhang,SPDSCombettes,SPDSChambolle}.
Roughly speaking, at each iteration, the algorithms activate only the operations associated with randomly chosen variables, so that the said costs are significantly reduced.
Actually, several studies show the utility of such block-coodinate-wise randomization in image restoration \cite{ICASSP2016SDCAADMM,SPDSCombettesOnlineIR,SPDSChambolle}. %

Incidentally, the above studies aim at \emph{unconstrained formulation}, i.e., minimizing the sum of a regularization and a data-fidelity term.
On the other hand, \emph{constrained formulation}, i.e., minimizing a regularization term subject to a hard constraint on data-fidelity has an important advantage over the unconstrained one in terms of facilitating the selection of regularization terms and hyperparameters, as addressed in \cite{ConstraintMI1982,CombettesConst1994,CSALSA,ConstStuck2011,ConstPoisson2012,IdivCons2013,EPIpre,TSP2015Involved}.
However, such a data-fidelity constraint is \emph{not separable} as is, i.e., it cannot be decomposed into block-coordinate-wise constraints, so that it does not suit randomized activation.

In this paper, we bridge the gap between the randomized nature of stochastic proximal methods and the non-separability of data-fidelity constraints by leveraging \emph{epigraphical projection} \cite{EPIpre,ICASSP2014}.
We focus on the $\ell_2$ data-fidelity constraint, and introduce its equivalent expression via certain epigraphs, which enables us to deal with the constrained formulation by stochastic proximal splitting algorithms with \emph{randomized} epigraphical projection.
Specifically, we develop a randomized solver for the problem based on a stochastic primal-dual hybrid gradient algorithm \cite{SPDSChambolle}.
The efficiency of our method is demonstrated on CT image reconstruction. 

\section{Preliminaries}
\subsection{Proximal Tools}
The \emph{proximity operator} \cite{Moreau} of index $\gamma>0$ of a proper lower semicontiuous convex function $f\in\Gamma_0(\R^N)$\footnote{
	The set of all proper lower semicontinuous convex functions on $\R^N$ is denoted by $\Gamma_0(\R^N)$.
} is defined as
\begin{equation*}
\prox_{\gamma f}:\R^N\rightarrow\R^N:\x\mapsto\argmin_{\y}f(\y)+\tfrac{1}{2\gamma}\|\y-\x\|^2.
\end{equation*}

The indicator function of a nonempty closed convex set $C$, denoted by $\iota_C$, is defined as 
\begin{equation*}
\iota_C(\x):=\begin{cases}
0 & \mbox{if } \x\in C\\
\infty & \mbox{otherwise.} 
\end{cases}
\end{equation*}
Since the function returns $\infty$ when the input vector is outside of $C$, it acts exactly as the hard constraint represented by $C$ in minimization.
The proximity operator of $\iota_C$ is equivalent to the (metric) projection onto $C$, i.e.,
\begin{equation*}
\prox_{\gamma\iota_C}(\y) = P_C(\y) := \argmin_{\x\in C}\|\y-\x\|.
\end{equation*}

\subsection{Stochastic Primal-Dual Hybrid Gradient Algorithm}\label{sec:SPDHG}
A stochastic primal-dual hybrid gradient algorithm (SPDHG) \cite{SPDSChambolle} was proposed to optimize the following problem:
\begin{equation}\label{prob:SPDHG}
\min_{\x\in\R^N} \textstyle\sum_{i=1}^I f_i(\A_i\x) + g(\x),
\end{equation}
where $f_i$ and $g$ are proper lower semicontinuous convex functions, and $\A_i$ are bounded linear operators (matrices).
Here we assume that the proximity operators of $f_i$ and $g$ are easy to compute.

SPDHG does not solve \eqref{prob:SPDHG} directly but solves the saddle point problem reformulated from \eqref{prob:SPDHG}, given by
\begin{equation}\label{prob:SPDHGsadd}
\min_{\x\in\R^N}\sup_{\y_i\in\R^{M_i}} \textstyle\sum_{i=1}^I \langle \A_i\x, \y_i \rangle - f_i^*(\y_i) + g(\x),
\end{equation}
where $f_i^*$ are the convex conjugate of $f_i$.
We note that the proximity operator of $f_i^*$ can be computed via that of $f_i$ as
$\prox_{\gamma f_i^*}(\x)=\x - \gamma\prox_{\frac{1}{\gamma}f_i}(\textstyle\frac{1}{\gamma}\x)$ \cite[Theorem~14.3(ii)]{Combettesbook}.

Let $\Sfrak\subset\{1,\ldots,I\}$ be a random subset of the indices of the dual variables $\y_i$ in \eqref{prob:SPDHGsadd},
and define $\A:=(\A_1,\ldots,\A_I)$, $\y:=(\y_1,\ldots,\y_I)$ and $\LAMBDA:=\mbox{diag}(p_1^{-1}\I,\ldots,p_I^{-1}\I)$ with the parameters $p_i$ being probabilities that an index is selected in each iteration.
Then SPDHG is formalized as follows: for given $\x^{0}$, $\y^{0}$, $\tau>0$, $\rho_i>0$, and $\overline\y^{(0)}=\y^{0}$, iterate
\begin{equation}\label{SPDS}
\left\lfloor\begin{array}{l}
\x^{(k+1)}=\prox_{\tau g}(\x^{(k)}-\tau\A^*\overline\y^{(k)}),\\
\mbox{Select } \Sfrak^{(k+1)}\subset\{1,\ldots,I\}\\
\y_i^{(k+1)}=\begin{cases}
\prox_{\rho_i f_i^*}(\y_i^{(k)} + \rho_i\A_i\x^{(k+1)}) & \mbox{if } i\in\Sfrak^{(k+1)}\\
\y_i^{(k)} & \mbox{otherwise}
\end{cases}\\
\overline{\y}^{(k+1)} = \y^{(k+1)} + \LAMBDA(\y^{(k+1)}-\y^{(k)})
\end{array}\right.
\end{equation}
We should note that the matrix-vector multiplication $\A^*\overline\y^{(k)}$ in \eqref{SPDS} can be computed by using only the selected $\A_i$ and the previous dual variable (see \cite[Remark 1 and 2]{SPDSChambolle} for details).
This means that each iteration requires both $\A_i$ and $\A_i^*$ to be evaluated only for each selected index $i\in\Sfrak^{(k+1)}$. 
With a mild condition on the stepsizes $\tau$ and $\rho_i$, the algorithm converges to an optimal solution of \eqref{prob:SPDHGsadd} almost surely in the sense of the Bregman distance (see \cite[Theorem~4.3]{SPDSChambolle} for details).

\section{Proposed Method}
\subsection{Problem Formulation}
Consider the following data observation model:
\begin{equation}\label{ir}
\v = \PHI\bar\u+\n,
\end{equation}
where $\bar\u\in\R^N$ is a latent signal we wish to estimate,
$\PHI\in\R^{M\times N}$ represents an observation process,
$\n\in\R^M$ is an additive white Gaussian noise,
and $\v\in\R^M$ is observed data.

Based on the model in \eqref{ir}, we aim at the following form of constrained signal reconstruction:
\begin{equation}\label{vir}
\min_{\u\in\R^N} \textstyle\sum_{j=1}^J\Rmath_j(\PSI_j\u) \mbox{ s.t. }\begin{cases}
\|\PHI\u-\v\|^2\leq\bar\varepsilon,\\
\u\in[\underline\mu,\overline\mu]^N,
\end{cases}
\end{equation}
where $\Rmath_j(\PSI_j\cdot)$ are regularization terms with functions $\Rmath_j\in\Gamma_0(\R^{P_j})$ and matrices $\PSI_j\in\R^{P_j\times N}$,
the first hard-constraint is $\ell_2$ delity with the radius $\bar\varepsilon>0$, and the second one is a range constraint.
We assume that the proximity operators of $\Rmath_j$ are available.

\subsection{Reformulation via Epigraphs}
Since the $\ell_2$ data-fidelity constraint in \eqref{vir} is not separable, we cannot directly solve the problem by stochastic algorithms with block-coordinate-wise randomization.
To circumvent this, we give another expression of the constraint as follows:
\begin{equation}\label{anoex}
\|\PHI\u-\v\|^2\leq\bar\varepsilon \Leftrightarrow \begin{cases}
\|\PHI_1\u-\v_1\|^2\leq\varepsilon_1,\\
\hspace{8mm}\vdots\\
\|\PHI_L\u-\v_L\|^2\leq\varepsilon_L,\\
\sum_{l=1}^L \varepsilon_l\leq\bar\varepsilon,
\end{cases}
\end{equation}
where $\varepsilon_l\in\R$ are additional variables, $(\PHI_1^\top\cdots\PHI_L^\top)^\top = \PHI$ and $(\v_1^\top\cdots\v_L^\top)^\top = \v$.
Let us define the epigraphs of $\v_l$-centered squared distance, denoted by $\Scal_l$, and a half space $\Vcal$, as 
\begin{align}
\mathcal{S}_l &:= \{(\x,\eta)\in\R^{Q_l}\times\R|\|\x-\v_l\|^2\leq\eta\}\label{epi}\\
\mathcal{V} &:= \{(\eta_1,\ldots,\eta_L)\in\R^L|\textstyle\sum_{l=1}^L\eta_l\leq\bar\varepsilon\},\label{half}
\end{align}
Then, with \eqref{anoex}, \eqref{epi}, and \eqref{half}, Problem~\eqref{vir} can be rewritten as
\begin{equation}\label{vir2}
\min_{\u,\varepsilon_1,\ldots,\varepsilon_I} \textstyle\sum_{j=1}^J\Rmath_j(\PSI_j\u) \mbox{ s.t. }\begin{cases}
(\PHI_1\u,\varepsilon_1)\in\mathcal{S}_1,\\
\hspace{8mm}\vdots\\
(\PHI_L\u,\varepsilon_L)\in\mathcal{S}_L,\\
\u\in[\underline\mu,\overline\mu]^N,\\
(\varepsilon_1,\ldots,\varepsilon_L)\in V.
\end{cases}
\end{equation}

By introducing the indicator functions of the hard constraints in \eqref{vir2}, which are denoted by $\iota_{\Scal_l}$, $\iota_{[\underline\mu,\overline\mu]}$, and $\iota_{\Vcal}$, we can further reformulate Problem~\eqref{vir2} as
\begin{align}
\min_{\u,\varepsilon_1,\ldots,\varepsilon_L} \textstyle\sum_{j=1}^J\Rmath_j(\PSI_j\u)+\textstyle\sum_{l=1}^L\iota_{\mathcal{S}_l}(\PHI_l\u,\varepsilon_l)\nonumber\\
+\iota_{[\underline\mu,\overline\mu]^N}(\u)+\iota_{\Vcal}(\varepsilon_1,\ldots,\varepsilon_L).\label{vir3}
\end{align}

\subsection{Algorithm}
Now we define $\boldsymbol{\varepsilon}:=(\varepsilon_1,\ldots,\varepsilon_L)$ and $\x:=(\u,\boldsymbol{\varepsilon})$, let $\e_1,\ldots,\e_L$ be the canonical basis of $\R^L$, and set $I:=J+L$,
\begin{align*}
&f_i(\A_i\x) := \Rmath_i(\PSI_i\u) \;\mbox{ for } i=1,\ldots,J,\\
&f_i(\A_i\x) := \iota_{\mathcal{S}_{i-J}}(\PHI_{i-J}\u, \e_{i-J}^\top\boldsymbol{\varepsilon})  \;\mbox{ for } i=J+1,\ldots,I,\\
&g(\x):=\iota_{[\underline\mu,\overline\mu]^N}(\u)+\iota_{\mathcal{V}}(\boldsymbol{\varepsilon}).
\end{align*}
Then \eqref{vir3} is reduced to \eqref{prob:SPDHG}, so that we can solve \eqref{vir3} by SPDHG.
We describe the whole algorithm in Algorithm~\ref{alg:RandEpi}.
\begin{remk}[Note on Algorithm~\ref{alg:RandEpi}]\label{remk1}\normalfont$\mbox{}$\\
(a) Our algorithm is designed to randomly pick up two indices respectively from the two index sets: one is $j\in\{1,\ldots,J\}$ with probability $1/J$, and the other is $l\in\{1,\ldots,L\}$ with probability $1/L$, so that both the $j$-th regularization term and the $l$-th data-fidelity epigraph are evaluated in each iteration.\\ 
(b) The matrix-vector multiplications required in each iteration are only for the selected indices $j$ and $l$.
This implies that the computational cost and memory requirement of Algorithm~\ref{alg:RandEpi} is much less than deterministic algorithms designed for solving \eqref{vir}.\\
(c) We give a stepsize setting rule as follows:
\begin{align*}
&\rho_\psi:=\textstyle\frac{\gamma}{\max_j\|\PSI_j\|}\;\;\forall j\in\{1,\ldots,J\},\\
&\rho_\phi:=\textstyle\frac{\gamma}{\max_l\|\PHI_l\|}\;\;\forall l\in\{1,\ldots,L\},\\
&\tau:=\textstyle\frac{\gamma}{\max\{J,L\}\max\{\max_j\|\PSI_j\|,\max_l\|\PHI_l\|\}},
\end{align*}
where $0<\gamma<1$. Similar stepsize setting is adopted in \cite{SPDSChambolle} but our rule is simpler.
\end{remk}

\begin{remk}[Projection computations in Algorithm~\ref{alg:RandEpi}]\normalfont$\mbox{}$\\
	Since the proximity operator of the indicator function of a nonempty closed convex set $C$ equals to the projection onto $C$, we need to compute $P_{[\underline\mu,\overline\mu]^N}$, $P_\Vcal$, and $P_{\Scal_i}$ in our algorithm.\\
	(a) The projection onto $[\underline\mu,\overline\mu]^N$ can be calculated by just pushing each entry of the input vector into $[\underline\mu,\overline\mu]$.\\
	(b) The projection onto $\mathcal{V}$ is given by \cite[(3.3-10)]{VSP}:
	\begin{equation*}
	P_{\mathcal{V}}(\bvarepsilon):=\begin{cases}
	\bvarepsilon & \mbox{if }\mathbf{1}^\top\bvarepsilon\leq\bar\varepsilon\\
	\bvarepsilon+\textstyle\frac{\bar\varepsilon-\mathbf{1}^\top\bvarepsilon}{L}\mathbf{1} & \mbox{otherwise},
	\end{cases}
	\end{equation*}
	where $\mathbf{1}$ is the all-one vector of size $L$.\\
	(c) The projection onto $\mathcal{S}_l$ is given as follows.
	\begin{prop}[Epigraphical projection of squared distance]
		Let $\z\in\R^N$ and let $\Scal := \{(\x,\eta)\in\R^{N}\times\R|\|\x-\z\|^2\leq\eta\}$.
		Then, for every $(\y,\zeta)\in\R^{N}\times\R$, by letting $d := \|\y-\z\|$, the projection onto $\Scal$ is given by
		\begin{equation}\label{epipro}
		P_\Scal(\y,\zeta) = (\alpha\y + (1-\alpha)\z, \max\{\alpha^2 d^2,\zeta\}),
		\end{equation}
		where
		\begin{align}
		&\hspace{-3mm}\alpha = \begin{cases}
		1 & d^2\leq\zeta,\\
		\frac{\beta}{d} & \mbox{otherwise,}
		\end{cases}\\
		&\hspace{-3mm}\textstyle\beta=(\frac{d}{4} + (\frac{d^2}{16} - (\frac{\zeta}{3} - \frac{1}{6})^3)^{\frac{1}{2}})^{\frac{1}{3}} + \frac{(\frac{\zeta}{3} - \frac{1}{6})}{(\frac{d}{4} + (\frac{d^2}{16} - (\frac{\zeta}{3} - \frac{1}{6})^3)^{\frac{1}{2}})^{\frac{1}{3}}}.\hspace{-2mm}\label{beta}
		\end{align}		
	\end{prop}
    \hspace{-5.3mm}Proof sketch: Equation~\eqref{epipro} can be easily obtained from \cite[Proposition 4]{EPIpre}.
    Using the same proposition, we can see that
    \begin{equation}\label{prf1}
    \alpha = d^{-1}\prox_{\frac{1}{2}(\max\{|\cdot|^2-\zeta,0\})^2}(d).
    \end{equation}
    Clearly, $\alpha=1$ when $d^2<\zeta$.
    When $d^2>\zeta$, by simple calculation, the solution of the proximity operator in \eqref{prf1} is reduced to the real solution of the following cubic equation:
    \begin{equation}\label{prf2}
    	2x^3 + (1-2\zeta)x - d = 0.
    \end{equation}
    Finally, applying Cardano formula to \eqref{prf2} yields \eqref{beta}.\footnote{One of the reviewers pointed out that the above result can also be proven by a combination of \cite[Proposition~5.1]{PesquetProxDiv} and \cite[Example~3.8]{CombettesPers}.}
\end{remk}\vspace{-4mm}
\begin{algorithm}[h]
	\LinesNumbered
	\SetKwInOut{Input}{input}
	\SetKwInOut{Output}{output}
	\SetKwInOut{Initialize}{initialize}
	\label{alg:RandEpi}
	\caption{Proposed algorithm for solving \eqref{vir}}
	\Input{$\u^{(0)},\bvarepsilon^{(0)},\z^{(0)},\w^{(0)},\t^{(0)},\bzeta^{(0)},\xii^{(0)}$}
	\Initialize{$\overline\t^{(0)}=\t^{(0)},\overline\xii^{(0)}=\xii^{(0)}$}
	\For{$k=0,\ldots,K-1$}{
		$\u^{(k+1)}=P_{[\underline{\mu,\overline{\mu}}]^N}(\u^{(k)}-\tau\overline\t^{(k)})$\;
		$\boldsymbol{\varepsilon}^{(k+1)}=P_{\mathcal{V}}(\boldsymbol{\varepsilon}^{(k)}-\tau\overline\xii^{(k)})$\;
		Select $j\in\{1,\ldots,J\}$ and $l\in\{1,\ldots,L\}$.\;
			$\tilde\z_j^{(k)}=\z_j^{(k)}+\rho_\psi\PSI_j\u^{(k+1)}$\;
			$\z_j^{(k+1)}=\tilde\z_j^{(k)} - \rho_\psi\prox_{\Rmath_j/\rho_\psi}(\tilde\z_j^{(k)}/\rho_\psi)$\;
			$\hat\z_j^{(k)}=\PSI_j^\top(\z_j^{(k+1)}-\z_j^{(k)})$\;
			$\tilde\w_l^{(k)}=\w_l^{(k)}+\rho_\phi\PHI_l\u^{(k+1)}$\;
			$\tilde\zeta_l^{(k)}=\zeta_l^{(k)}+\rho_\phi\varepsilon_l^{(k+1)}$\;
			$(\w_l^{(k+1)},\zeta_l^{(k+1)})=(\tilde\w_l^{(k)},\tilde\zeta_l^{(k)}) - \rho_\phi P_{\mathcal{S}_l}((\tilde\w_l^{(k)},\tilde\zeta_l^{(k)})/\rho_\phi)$\;
			$\hat\w_l^{(k)}=\PHI_l^\top(\w_l^{(k+1)}-\w_l^{(k)})$\;
			$\hat\bzeta_l^{(k)}=\e_l(\zeta_l^{(k+1)}-\zeta_l^{(k)})$\;
		$\t^{(k+1)} = \t^{(k)} + \hat\z_j^{(k)} + \hat\w_l^{(k)}$\;
		$\overline\t^{(k+1)} = \t^{(k+1)} + (1+J)\hat{\z_j}^{(k)} + (1+L)\hat{\w_l}^{(k)}$\;
		$\xii^{(k+1)} = \xii^{(k)} + \hat\bzeta_l^{(k)}$\;
		$\overline\xii^{(k+1)} = \xii^{(k+1)} + (1+L)\hat{\bzeta_l}^{(k)}$
	}
\end{algorithm}\vspace{-5mm}

\begin{figure*}[t]
	\begin{minipage}{0.33\hsize}
		\centering\includegraphics[width=\hsize]{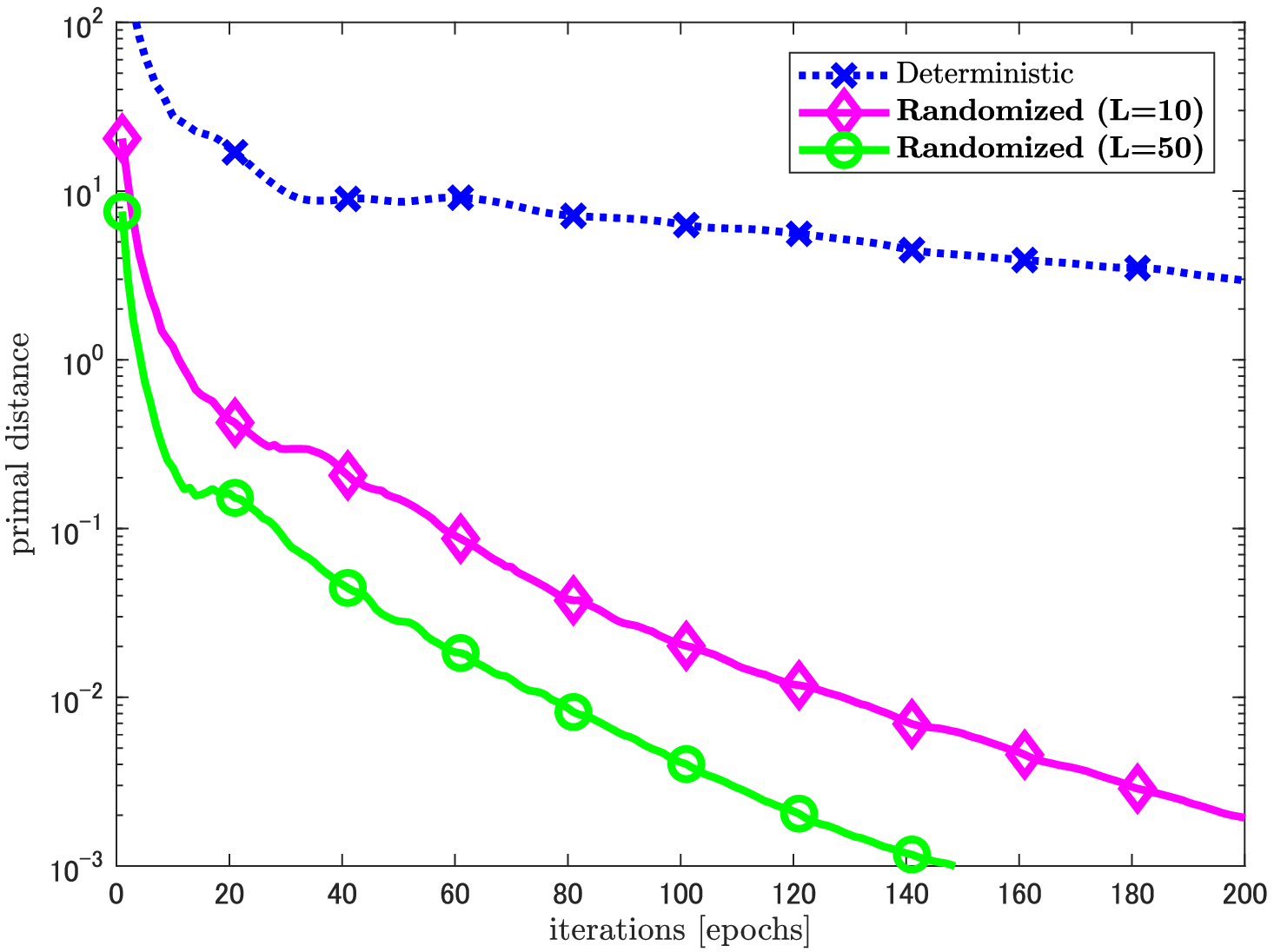}
	\end{minipage}
	\begin{minipage}{0.33\hsize}
		\centering\includegraphics[width=\hsize]{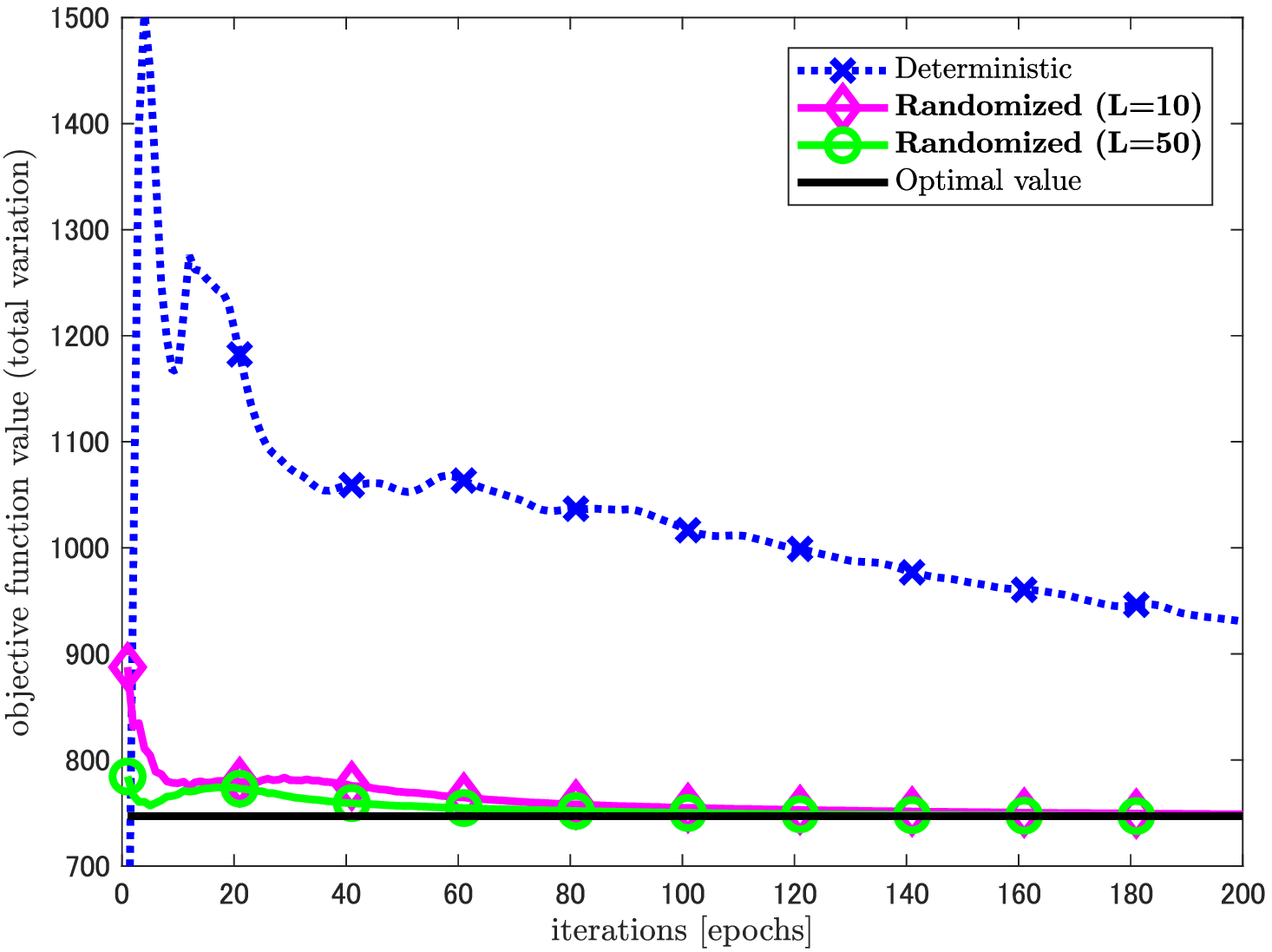}
	\end{minipage}
	\begin{minipage}{0.33\hsize}
		\centering\includegraphics[width=\hsize]{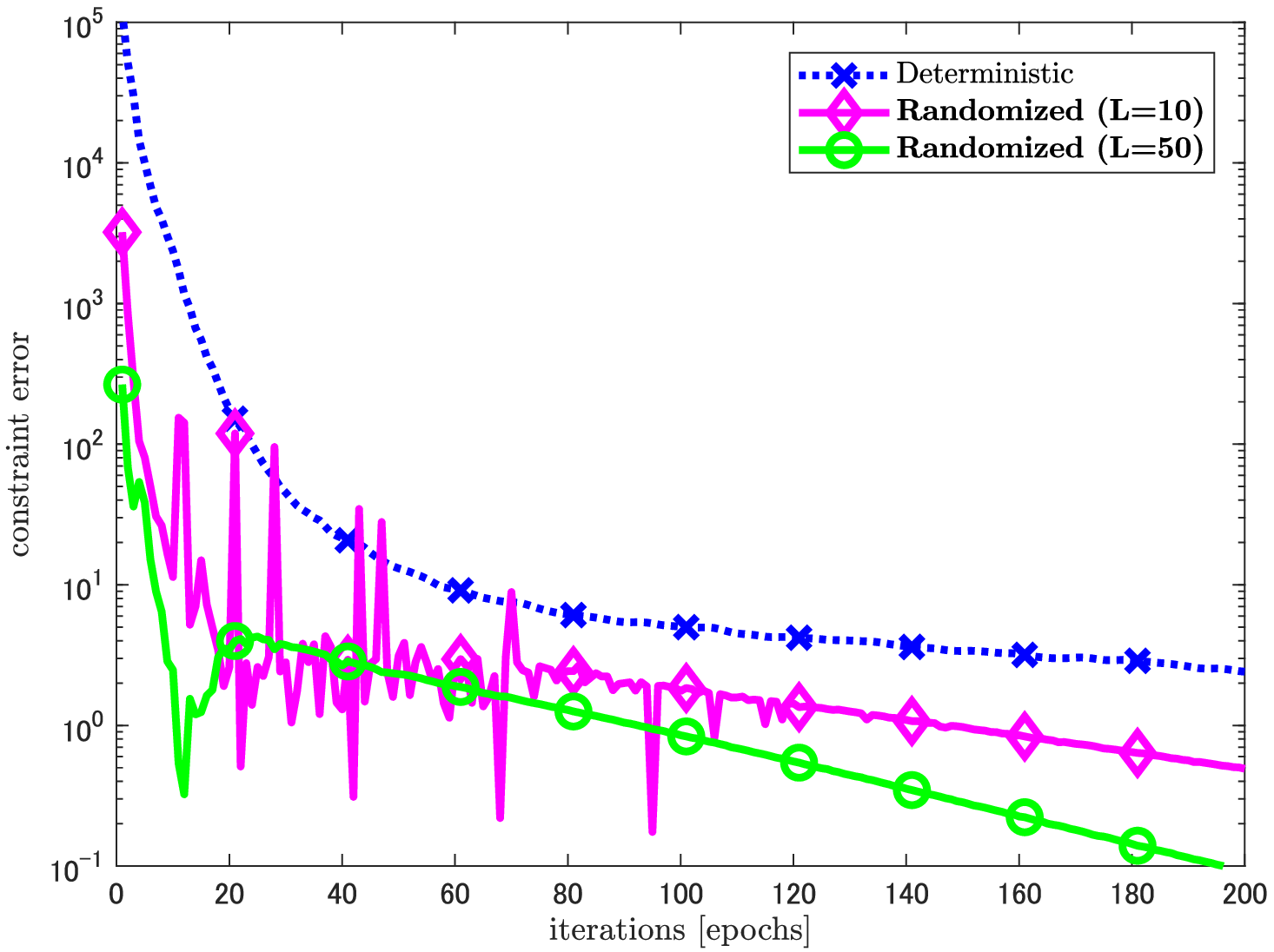}
	\end{minipage}
	\vspace{-4mm}
	\caption{\small{Convergence profiles of Algorithm~\ref{alg:RandEpi} (``Randomized'') and its deterministic counterpart (``Deterministic'') on CT image reconstruction in terms of the primal distance (left), objective function value (center), and constraint error (right). Note that the optimal value of the objective function (the black line in the center figure) was measured on $\u^\star$.}}
	\label{fig:CT}
	\vspace{-2mm}
\end{figure*}
\begin{figure*}[t]
	\begin{minipage}{0.196\hsize}
		\centering\includegraphics[width=\hsize]{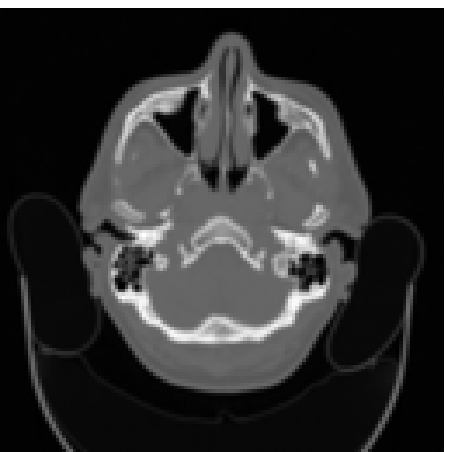}
		\centerline{\footnotesize{Original ($\bar{\u}$)}}\vspace{-1mm}
		\centerline{\footnotesize{$\;$}}
	\end{minipage}
	\begin{minipage}{0.196\hsize}
		\centering\includegraphics[width=\hsize]{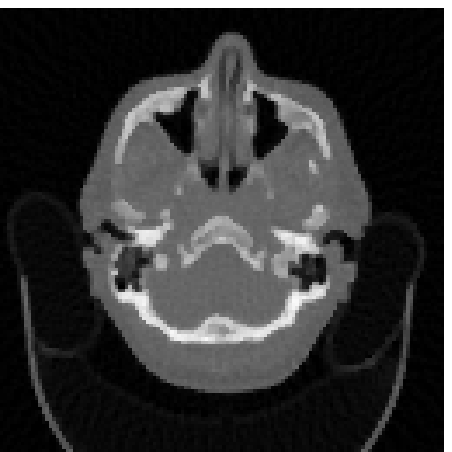}
		\centerline{\footnotesize{Deterministic}}\vspace{-1mm}
		\centerline{\footnotesize{PSNR=34.22 [dB]}}
	\end{minipage}
	\begin{minipage}{0.196\hsize}
		\centering\includegraphics[width=\hsize]{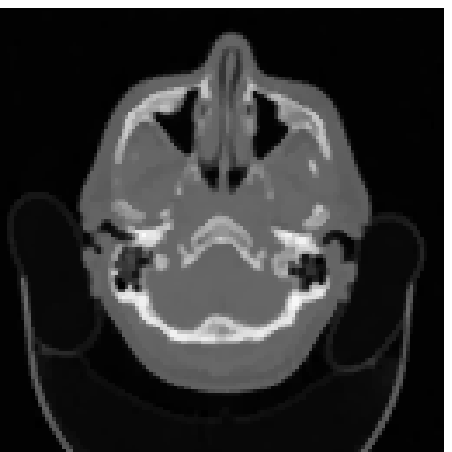}
		\centerline{\footnotesize{\textbf{Stochastic} ($L=10$)}}\vspace{-1mm}
		\centerline{\footnotesize{PSNR=37.51 [dB]}}
	\end{minipage}
	\begin{minipage}{0.196\hsize}
		\centering\includegraphics[width=\hsize]{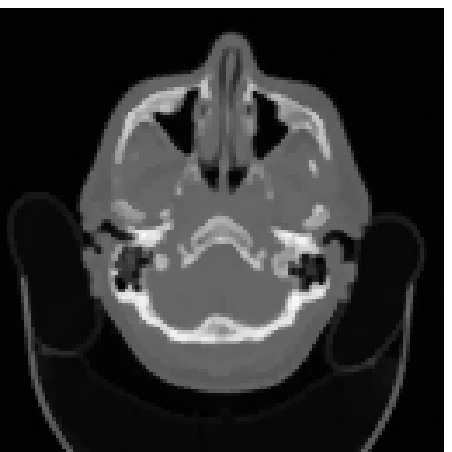}
		\centerline{\footnotesize{\textbf{Stochastic} ($L=50$)}}\vspace{-1mm}
		\centerline{\footnotesize{PSNR=37.47 [dB]}}
	\end{minipage}
	\begin{minipage}{0.196\hsize}
	\centering\includegraphics[width=\hsize]{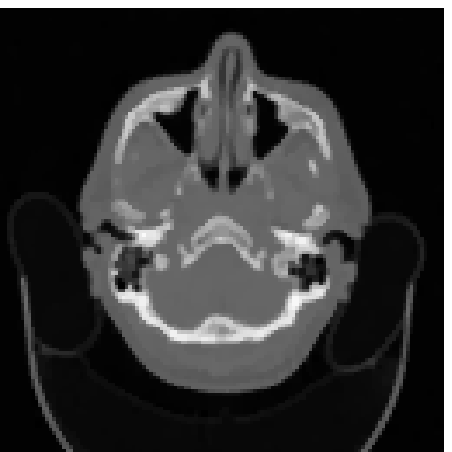}
	\centerline{\footnotesize{Optimal ($\u^\star$)}}\vspace{-1mm}
	\centerline{\footnotesize{PSNR=37.47 [dB]}}
\end{minipage}
	\vspace{-4mm}
	\caption{\small{Resulting images on CT image reconstruction (200 iterations [epochs]).}}
	\label{fig:CTresult}
	\vspace{-5mm}
\end{figure*}

\section{Numerical Experiments}
We examined the performance of the proposed method by comparing it with the deterministic counterpart, i.e., the deterministic primal-dual hybrid gradient algorithm \cite{PDChambolle} on CT image reconstruction.
All experiments were performed using MATLAB (R2017b), on a Windows 10 Pro laptop computer with an Intel Core i7 2.1 GHz processor and 16 GB of RAM.

For the original image $\bar{\u}$ in \eqref{ir}, we used a head CT scan image of size $128\times128$ ($N=16384$) picked up from the CT dataset \cite{CTdata}.
The matrix $\PHI$ in \eqref{ir} was set to a parallel beam projection (Radon transform) matrix with 60 projection angles.
We would like to note that the nonzero entries of $\PHI$ account for about $1.4\%$ of all the entries, i.e., $\PHI$ is sparse,
so that this is advantageous for the deterministic algorithm, compared with the cases of dense $\PHI$.
We also note that we use a small image because the deterministic algorithm has to load full $\PHI$ of size $M\times N$ ($M=11100, N=128^2$) in each iteration.
The observed data was generated by adding white Gaussian noise with standard deviation $\sigma=10/255$ to $\PHI\bar{\u}$.

A full CT image was estimated by constrained total variation (TV) minimization, which is a special case of Problem~\eqref{vir}.
Specifically, we employed the anisotropic TV \cite{ROF} for the regularization function $\Rmath_j\circ\PSI_j$ in Prob.~(\ref{vir}).
In this case, $J = 2$, and the matrix $\PSI_1$ and $\PSI_2$ are equal to the vertical and horizontal discrete gradient operators $\D_v$ and $\D_h$ with Neumann boundary, respectively.
Both $\Rmath_1$ and $\Rmath_2$ are the $\ell_1$ norm, and its proximity operator can be calculated by a simple $\mathcal{O}(N)$ soft-thresholding operation.
We adopted an eight-bit dynamic range constraint $[0,255]^N$ for $[\underline\mu,\overline\mu]^N$.
For a fair comparison, the parameter $\bar{\varepsilon}$ was set to an oracle value, i.e., $\|\PHI\bar{\u}-\n\|$.
For the stepsizes of the deterministic algorithm, we employed the setting rule suggested in \cite{PDChambolle}.
For our algorithm, see Remark~\ref{remk1}(c) (Note: we choose $\gamma=0.99$).


We adopted the following three convergence criteria:\\
(i) Primal distance: the squared distance between the current estimate $\u^{(k)}$ and an optimal solution $\u^\star$, i.e., $\|\u^{(n)} - \u^\star\|^2$. Since $\u^\star$ is analytically unavailable, it was pre-computed by the deterministic algorithm with $2\times10^5$ iterations.\\
(ii) Objective function value: the value of TV defined by $\|\D_v\u^{(k)}\|_1+\|\D_h\u^{(k)}\|_1$.\\
(iii) Constraint error: the absolute difference between $\bar\varepsilon$ and $\|\PHI\u^{(k)}-\v\|^2$. Since any optimal solution $\u^\star$ satisfies $\|\PHI\u^\star-\v\|^2=\bar\varepsilon$, this value should converge to zero.

The left of Fig.~\ref{fig:CT} shows the convergence of the primal distance,
where "iterations [epochs]" means that the number of iterations is divided by $L$. 
We examined the cases of $L=10,50$ in this experiment.
One sees that the proposed method (``Randomized'') converges much faster than the deterministic counterpart (``Deterministic'').
Similar convergence behavior can be observed in the center and right of Fig.~\ref{fig:CT},
where the convergence profiles of the objective function value and the constraint error are plotted, respectively.
The resulting images are depicted in Fig.~\ref{fig:CTresult}, which illustrates that our algorithm properly works.\vspace{-3mm}

\section{Conclusion}\vspace{-3mm}
We have proposed an efficient constrained signal reconstruction framework based on a stochastic primal-dual splitting algorithm with randomized epigraphical projection.
Since the proposed method does not require the multiplication of full $\PHI$ and variables in each iteration,
it would be a powerful choice when $\PHI$ is large and not structured while keeping the benefits of the constrained formulation.

In this paper, we discuss only the $\ell_2$ data-fidelity case but this framework can also be applied to other data fidelity constraints, for example, the $\ell_1$ case, as long as their epigraphical projections are computable. 
Also, with a slight extension, our method would be able to handle signal decomposition models, such as image decomposition \cite{CTdecAujol2006,TIP2014}.

\balance
\footnotesize{
\bibliographystyle{IEEEbib}
\bibliography{ICASSP2019_RandEpi.bib}
}

\end{document}